\input amstex
\documentstyle{amsppt}
\magnification=1100

\def\aa{\alpha}
\def\ee{\epsilon}
\def\ll{\lambda}
\def\oo{\omega}
\def\BB{\Cal B}
\def\CC{\Cal C}
\def\DD{\Cal D}
\def\HH{\Cal H}
\def\LL{\Cal L}
\def\ZZ{\Bbb Z^n}
\def\RR{\Bbb R^n}

\topmatter

\title Covering lattice points by subspaces\endtitle
\author Imre B\'ar\'any, Gergely Harcos,
J\'anos Pach, G\'abor Tardos\endauthor

\address R\'enyi Institute of the Hungarian Academy of Sciences,
POB 127, H-1364 Budapest, Hungary and
Department of Mathematics, University College London,
Gower Street, London WC1E 6BT, England\endaddress
\email barany\@math-inst.hu\endemail
\thanks First author supported by grant T020914 of the
Hungarian National Foundation for Scientific Research (OTKA)\endthanks

\address Department of Mathematics, Princeton University, Fine Hall,
Washington Road, Princeton, NJ 08544, USA\endaddress \email
gharcos\@math.princeton.edu\endemail
\thanks Second author supported by grant 220/1762 of
Soros Foundation, Budapest\endthanks

\address R\'enyi Institute of the Hungarian Academy of Sciences,
POB 127, H-1364 Budapest, Hungary and
Courant Institute, 251 Mercer Street, New York, NY 10012, USA\endaddress
\email pach\@cims.nyu.edu \endemail
\thanks Third author supported by NSF grant CCR-9732101,
a PSC-CUNY Research Award, and grant T020914 of the Hungarian
National Foundation for Scientific Research (OTKA)\endthanks

\address R\'enyi Institute of the Hungarian Academy of Sciences,
POB 127, H-1364 Budapest, Hungary\endaddress
\email tardos\@math-inst.hu\endemail
\thanks Fourth author supported by grants T029255 and T030059
of the Hungarian National Foundation for Scientific Research (OTKA),
and grant FKFP 0607/1999 of the Hungarian Ministry of Education\endthanks


\abstract We find tight estimates for the
minimum number of proper subspaces needed to cover all lattice
points in an $n$-dimensional convex body $\CC$, symmetric
about the origin $0$. This enables us to prove the following
statement, which settles a problem of G. Hal\'asz. The maximum
number of $n$-wise linearly independent lattice points in the
$n$-dimensional ball $r\BB^n$ of radius $r$ around $0$ is
$O(r^{n/(n-1)})$. This bound cannot be improved.
We also show that the order of magnitude of the number
of different $(n-1)$-dimensional subspaces induced by the
lattice points in $r\BB^n$ is $r^{n(n-1)}$.
\endabstract

\endtopmatter

\document

\head 1. Introduction and statement of results\endhead

This paper was inspired by the following question of
G. Hal\'asz. {\it What is the maximal cardinality of a subset $S$ of
$r\BB^n\cap\ZZ$ such that all $n$-element subsets of $S$ are linearly
independent?} (Here $\BB^n$ denotes the unit ball around the
origin in $\RR$.) As any system of proper subspaces that cover
$r\BB^n\cap\ZZ$ provides an upper bound on the above quantity,
we would like to determine the size of the {\it smallest}
such covering system. We look
at these questions from a somewhat broader perspective.

We introduce the following notations. Let $\CC\subseteq\RR$ be
a convex compact body symmetric with respect to the origin.
For $1\le i\le n$, let $\ll_i$ denote the $i$-th {\it successive
minimum} of $\CC$. That is,
$$\ll_i=\min\{\ll|\dim(\ll\CC\cap\ZZ)\ge i\}.$$
Let $g(\CC)$ denote the minimum number of
proper subspaces covering $\CC\cap\ZZ$,
and let $h(\CC)$ denote the maximum number of points that
can be chosen from $\CC\cap\ZZ$ so that they are in
{\it general position}, i.e., no $n$ of them are linearly
dependent. Clearly, we have $h(\CC)\leq (n-1)g(\CC)$.

The following two theorems, providing a lower bound on
$h(\CC)$ and an upper bound on $g(\CC)$, respectively, give fairly
tight estimates for these quantities.

\proclaim{Theorem 1} If $\ll_n\leq1$ then
$$h(\CC)\ge{1-\ll_n\over16n^2}\min_{0<m<n}(\ll_m\dots\ll_n)
^{-\frac{1}{n-m}}.$$
\endproclaim

\proclaim{Theorem 2} If $\ll_n\leq 1$ then
$$g(\CC)\leq c2^nn^2\log n\min_{0<m<n}(\ll_m\dots\ll_n)^{-\frac{1}{n-m}},$$
where $c$ is an absolute constant.
\endproclaim

In Hal\'asz' question, $\CC$ is the $n$-dimensional ball,
$r\BB^n$, of radius $r>1$ around the origin, whose successive
minima satisfy $\ll_1=\ll_2=\ldots=\ll_n=1/r$. Thus, in this case,
Theorems 1 and 2 immediately imply that the correct orders of
magnitude of both $g(r\BB^n)$ and $h(r\BB^n)$ are $O(r^{n/(n-1)})$.

\remark{Remark 1} If $\ll_n>1$, then $g(\CC)=1$ and hence $h(\CC)<n$.
If $\ll_n<1-\ee$, by Theorems 1 and 2 the values of $g(\CC)$ and $h(\CC)$
are determined by the successive minima of $\CC$
up to a constant factor depending on $\ee$ and the
dimension $n$. For $\ll_n=1$ no such approximation is possible.
For arbitrary large $x>1$, consider the convex bodies
$$\CC_x=[-x,x]^{n-1}\times[-1,1]$$
and
$$\CC_x'=\hbox{conv}(\{-xe_i,xe_i|1\le i<n\}\cup\{-e_n,e_n\}),$$
where $(e_1,\ldots,e_n)$ is the standard basis of $\ZZ$. Both bodies have
the same sequence of successive minima: $\ll_i=1/x$ for $i<n$ and
$\ll_n=1$. However, $g(\CC_x)\ge2x$ and $h(\CC_x)\ge x/2,$
while $g(\CC_x')=2$ and $h(\CC_x')=n$.
\endremark

\remark{Remark 2} The integer lattice $\ZZ$ plays no particular role
in the above theorems. Our inequalities are preserved by affine
transformations, therefore they hold for $n$-dimensional lattices in general.
\endremark

\medskip

For any $r\geq 1$, let $\HH_r$ denote the set of all
$(n-1)$-dimensional subspaces (hyperplanes through $0$)
which contain $n-1$ linearly independent lattice
points from the ball of radius $r$ centered at the origin.

\proclaim{Theorem 3} There exist suitable positive constants $c_1$ and $c_2$,
depending only on $n$, such that
$$c_1r^{n(n-1)}\leq|\HH_r|\leq c_2r^{n(n-1)},$$
provided that $r$ is large enough.
\endproclaim

Let $r\geq 1$. Theorem 3 can be used to bound
$$s_r=\frac{1}{|\HH_r|}\sum_{H\in\HH_r}|H\cap r\BB^n\cap\ZZ|,$$
the average number of lattice points in $r\BB^n$ in
the hyperplanes belonging to $\HH_r$.

\proclaim{Corollary} 
The average $s_r$ is bounded by a constant depending on the dimension $n$.
\endproclaim

\remark{Remark 3} By analyzing the dependence of $c_1$ on $n$ it is possible
to show that $s_r\leq 2^{n^3+O(n^2\log n)}$.
\endremark

\medskip

In Section 2, we essentially show that within $\CC\cap\ZZ$
one can represent a finite projective space over a
relatively small prime (see Lemma). To establish Theorem 1,
we combine this result with a well known construction of
P. Erd\H os (see \cite{11, Appendix}).

Section 3 contains the proof of Theorem 2.
This proof is also constructive: in most cases, to cover
$\CC\cap\ZZ$ we take all subspaces perpendicular to an
integer vector in a body homothetic to the polar of $\CC$.

The proofs of Theorem 3 and the Corollary are given in Section 4.

\bigskip

The related (but different) problem of covering the lattice points
within a convex body by {\it affine} subspaces was first
investigated by K. Bezdek and T. Hausel \cite{2}.
They only considered
1-codimensional subspaces, i.e. hyperplanes (as we do here).
Their work was sharpened and extended to the general case
by I. Talata \cite{14}. The estimates in these two papers
are given in terms of the dimension $n$ and the lattice width
of the convex body.

\head 2. Proof of Theorem 1\endhead

The proof is based on the following

\proclaim{Lemma} Let
$\ll_n<1$ and suppose that $p$ is an integer satisfying
$$1<p<{1-\ll_n\over8n^2}\min_{0<m<n}(\ll_m\dots\ll_n)^{-\frac{1}{n-m}}.$$

Then, for any $v\in\RR$, there exist an integer
$1\le j<p$ and a lattice point $w\in\ZZ$ with $jv+pw\in\CC$.
\endproclaim

\demo{Proof of Lemma}
Find linearly independent vectors $v_i\in\ll_i\CC\cap\ZZ$
for $i=1,\ldots,n$. Any vector $x\in\RR$ can be uniquely written in the
form $x=\sum_{i=1}^na_iv_i+\sum_{i=1}^nb_iv_i$ with $a_i\in\Bbb Z$ and
$b_i\in(-1/2,1/2]$. Here $\sum_{i=1}^na_iv_i\in\ZZ$ and
$$\sum_{i=1}^nb_iv_i\in
\hbox{conv}\left\{{v_i\over2p\ll_i},-{v_i\over2p\ll_i}
\Big|1\le i\le n\right\}\subseteq{\CC\over2p},$$
whenever $\sum_{i=1}^n\ll_i|b_i|\le1/(2p)$.
Thus, the density $d$ of the periodic set
$$S=\frac{\CC}{2p}+\ZZ$$
is at least the probability that for independent uniform random
numbers $b_i\in[0,1/2]$ we have $\sum_{i=1}^n\ll_ib_i\le1/(2p)$. This
inequality is satisfied if $\ll_ib_i\le\ee/(2pn)$ for all $i<n$ and
$\ll_nb_n<(1-\ee)/(2p)$, where $\ee=(1-\ll_n)/2$.
Thus, we have
$$d\ge\min\left(1,\frac{1-\ee}{p\ll_n}\right)
\prod_{i=1}^{n-1}\min\left(1,\frac{\ee}{pn\ll_i}\right).$$
This lower bound on $d$ takes the form
$$A_m=\prod_{m\le i<n}{\ee\over pn\ll_i}$$
or
$$B_m={1-\ee\over p\ll_n}\prod_{m\le i<n}{\ee\over pn\ll_i} \ ,$$
where $1\le m\le n$ is an appropriate integer (the product is empty in
case $m=n$).

We claim that each of these values is larger than $1/p$, so we
have $d>1/p$. The inequality $B_m>1/p$ is equivalent to
$$p^{n-m}<{(1-\ee)\ee^{n-m}\over n^{n-m}\ll_m\ldots\ll_n}.$$
This is true, by the choice of $\ee$, for $m=n$, and, by our bound on
$p$, otherwise. The inequality $A_m>1/p$ is equivalent to
$$C_m=p^{n-m-1}\ll_m\ldots\ll_{n-1}<\left({\ee\over n}\right)^{n-m}.$$
If $m=n$, this is true, because $p>1$. Suppose $m<n$,
and use our bound on $p$ to get
$$C_m<\frac{1}{p\ll_n}\left(\frac{\ee}{4n^2}\right)^{n-m}.$$
If $\ll_n\geq 1/2$ then $p\ll_n\geq 1$, hence the desired inequality follows.
If $\ll_n<1/2$ then $\ee>1/4$, hence the previous inequality yields
$$C_m<\frac{1}{p\ll_n}\left(\frac{\ee}{n}\right)^{n-m+1}.$$
On the other hand,
using the monotonicity of the sequence $(\ll_i)$, we obtain
$$C_m\le p^{n-m-1}\ll_n^{n-m}<(p\ll_n)^{n-m}.$$
Taking a weighted
geometric mean of the last two bounds, we get
$$C_m<\left\{\frac{1}{p\ll_n}\left(\frac{\ee}{n}\right)^{n-m+1}\right\}
^\frac{n-m}{n-m+1}\left\{(p\ll_n)^{n-m}\right\}^\frac{1}{n-m+1}=
\left(\frac{\ee}{n}\right)^{n-m},$$
as required. This proves $A_m>1/p$ and hence $d>1/p$.

Consider the periodic sets $S+jv/p$ for $j=0,\ldots,p-1$. Each of these $p$
sets have density $d>1/p$ thus two of these sets must intersect. We have
$$\frac{j_1v}{p}+\frac{u_1}{2p}+w_1=\frac{j_2v}{p}+\frac{u_2}{2p}+w_2,$$
for some $0\le j_1<j_2<p$, some
$u_1,u_2\in\CC$ and some $w_1,w_2\in\ZZ$. For $1\le j=j_2-j_1<p$ and
$w=w_2-w_1\in\ZZ$, we have
$$jv+pw=\frac{u_1-u_2}{2}\in\CC,$$
verifying the statement of the Lemma. \qed
\enddemo

Now it is easy to finish the proof of Theorem 1.
Let $p$ be the largest prime number satisfying the condition
in the Lemma. If such a prime does not exist,
then the statement of the theorem is trivial. The points of the
{\it discrete moment curve} (used by Erd\H os in connection
with Heilbronn's triangle problem \cite{11}),
$v_i=(1,i,i^2,\ldots,i^{n-1})$ for integer values
$0\le i<p$ (and $v_\infty=(0,\ldots,0,1)\in\ZZ$) are $n$-wise linearly
independent over the $p$-element field.
By the Lemma, we have integers $1\le j_i<p$ and integer vectors $w_i$ with
$v_i'=j_iv_i+pw_i\in\CC$. Clearly, the vectors $v_i'$ are integer
vectors, and they are $n$-wise linearly indepent over the $p$-element field,
and hence over the reals. This shows $h(\CC)>p$, and an application of
Chebyshev's theorem on prime numbers concludes the proof.

\head 3. Proof of Theorem 2\endhead

Let $\CC^0$ denote the {\it polar body} of $\CC$, i.e.,
$$\CC^0=\{x\in\RR:ux\leq 1\text{ for all $u\in\CC$}\}.$$
Denote by
$\mu_1\leq\dots\leq\mu_n$ the successive minima
of $\CC^0$. It is known that
$$1\leq\ll_i\mu_{n-i+1}\leq c_1n\log n\quad(i=1,\dots,n)$$
where $c_1$ is an absolute constant. The lower bound is a classical inequality
of Mahler \cite{10}, the upper one has been recently proved by
Banaszczyk \cite{1}.

Fix any integer $0<m<n,$ for the rest of the argument. It follows that
$$1\leq(\ll_m\dots\ll_n)(\mu_1\dots\mu_{n-m+1})\leq(c_1n\log n)^{n-m+1}.\tag 1$$
For technical reasons, we will consider any increasing sequence
$$0<\nu_1<\dots<\nu_{n-m+1}$$
such that no ratio $\nu_i/\nu_j$ ($i\neq j$) is rational and
$$\mu_i\leq\nu_i\quad(i=1,\dots,n-m+1).$$
Let
$$w_i\in\mu_i\CC^0\cap\ZZ\quad(i=1,\dots,n-m+1)$$
be linearly independent vectors, and consider some sets of
integer vectors of the form
$$\DD_\aa^+=\left\{\sum_{i=1}^{n-m+1}a_iw_i:a_i\in[0,\aa/\nu_i]\cap\Bbb Z\right\},$$
$$\DD_\aa=\left\{\sum_{i=1}^{n-m+1}a_iw_i:a_i\in[-\aa/\nu_i,\aa/\nu_i]\cap\Bbb Z\right\},$$
where $\aa$ is a non-negative parameter to be specified later.
Clearly, $\DD_\aa$ is the union of $2^{n-m+1}$ isometric copies of $\DD_\aa^+$ satisfying
$$\DD_\aa\subseteq(n-m+1)\aa\CC^0\cap\ZZ.$$
Also, the difference of any two vectors from $\DD_\aa^+$ lies in $\DD_\aa$.
Let $f(\aa)$ be the number of points in the first set, i.e.,
$$f(\aa)=\left|\DD_\aa^+\right|=\prod_{i=1}^{n-m+1}\left(\left\lfloor\frac{\aa}{\nu_i}\right\rfloor+1\right).$$
Notice that $f(\aa)$ is an increasing, right continuous function which changes
by a factor of at most 2 at its points of discontinuity, i.e., for any $\aa>0$,
$$f(\aa)\leq 2f(\aa-).\tag 2$$
Also, $f(0)=1$ and
$$f(\aa)\geq\prod_{i=1}^{n-m+1}\frac{\aa}{\nu_i}.\tag 3$$
We claim that, whenever
$$f(\aa)>2(n-m+1)\aa+1\tag 4$$
holds, every lattice point in $\CC$ is perpendicular to some non-zero element of $\DD_\aa$.
To see this, fix any $u\in\CC\cap\ZZ$ and consider all the scalar products
$uv$ where $v\in\DD_\aa^+$. These scalar products are integers,
whose absolute values do not exceed $(n-m+1)\aa$.
Therefore, (4) implies the existence of two distinct
$v_1,v_2\in\DD_\aa^+$ with $uv_1=uv_2$.
Hence, the non-zero vector $v=v_1-v_2\in\DD_\aa$ is perpendicular to $u$.
We established that (4) implies
$$g(\CC)\leq |\DD_\aa|\leq 2^{n-m+1}f(\aa).\tag 5$$
By the right continuity of $f(\aa),$ there is a minimum $\aa$ such that
$$f(\aa)\geq 16(n-m+1)^\frac{n-m+1}{n-m}(\nu_1\dots\nu_{n-m+1})^\frac{1}{n-m}.$$
By (3), this $\aa$ satisfies
$$\aa\leq 4(n-m+1)^\frac{1}{n-m}(\nu_1\dots\nu_{n-m+1})^\frac{1}{n-m}.$$
In particular, we have
$$4(n-m+1)\aa\leq f(\aa).$$
The inequality $0<\ll_m\leq\dots\leq\ll_n\leq 1$ combined with (1) guarantees that
$$1\leq\mu_1\dots\mu_{n-m+1}\leq\nu_1\dots\nu_{n-m+1},$$
whence also
$$32\leq f(\aa).$$
The last two estimates on $f(\aa)$ show that (4) is satisfied.
In particular, $\aa>0$, therefore (5) combined with (2) yields
$$g(\CC)\leq 2^{n-m+2}f(\aa-)<
2^{n-m+6}(n-m+1)^\frac{n-m+1}{n-m}(\nu_1\dots\nu_{n-m+1})^\frac{1}{n-m}.$$
Taking the infimum of the right hand side over all admissible choices of the
sequence $0<\nu_1<\dots<\nu_{n-m+1}$, we get
$$\align g(\CC)
&\leq 2^{n-m+6}(n-m+1)^\frac{n-m+1}{n-m}(\mu_1\dots\mu_{n-m+1})^\frac{1}{n-m}\\
&\leq 2^{n-m+7}n(\mu_1\dots\mu_{n-m+1})^\frac{1}{n-m}.\endalign$$
Combining this with (1), we obtain
$$\align g(\CC)
&\leq 2^{n-m+7}n(c_1n\log n)^\frac{n-m+1}{n-m}(\ll_m\dots\ll_n)^\frac{1}{n-m}\\
&\leq 2^{n+7}c_1^2n^2\log n\bigl\{2^{-m}(n\log n)^\frac{1}{n-m}\bigr\}(\ll_m\dots\ll_n)^\frac{1}{n-m}.
\endalign$$
Here
$$2^{-m}(n\log n)^\frac{1}{n-m}\leq
\max\bigl\{(n\log n)^{2/n},2^{-n/2}n\log n\bigr\}$$
is bounded from above by an absolute constant, hence we can see that
$$g(\CC)\leq 2^ncn^2\log n(\ll_m\dots\ll_n)^\frac{1}{n-m},$$
where $c$ is some absolute constant. Minimizing over all integers $0<m<n,$
Theorem 2 follows.

\head 4. Proof of Theorem 3\endhead

The upper bound follows at once by noting that
$$|\HH_r|\leq{|r\BB^n\cap\ZZ|\choose n-1}={O(r^n)\choose n-1}=O(r^{n(n-1)}).$$

For any primitive integer vector $v$,
let $\LL(v)$ stand for the $(n-1)$-dimensional lattice
$\ZZ \cap v^{\perp}$ orthogonal to $v$, with
determinant $\det\LL(v)=|v|$. Write
$\ll_1(v)\leq\dots\leq\ll_{n-1}(v)$ for the successive minima
of $\LL(v)$, i.e.,
$$\ll_i(v)=\min\{\ll|\dim(\ll\BB^{n}\cap\LL(v))\ge i\}.$$
Denote by $\oo_n$ the volume of the unit ball $\BB^n$. According
to Minkowski's second fundamental theorem, we have
$$\ll_1(v)\dots\ll_{n-1}(v)\leq 2^{n-1}\oo_{n-1}^{-1}|v|.\tag 6$$

Define a set $V$ by
$$V=\{v \in \ZZ \: v \text{ is primitive and } |v| \leq \rho \},$$
where $\rho$ will be specified later.

\proclaim{Claim}
If $\rho$ is large enough, there are at least $\oo_n\rho^n/10$
elements $v \in V$ such that
$\ll_1(v) \geq D\rho^{\frac{1}{n-1}}$, where $D>0$ is
a suitable constant depending on $n$.
\endproclaim

Before proving the Claim, we show how it implies the lower bound in
Theorem 3. By (6),
whenever $\ll_1(v) \geq D\rho^{\frac{1}{n-1}},$ we have
$$
\ll_{n-1}(v) \leq 2^{n-1}\oo_{n-1}^{-1}|v|(D\rho^{\frac{1}{n-1}})^{-(n-2)}
\leq 2^{n-1}\oo_{n-1}^{-1}D^{-(n-2)}\rho^{\frac{1}{n-1}}.
$$
So, for at least $\oo_n\rho^n/10$ elements $v \in V$, $\LL(v)$
contains $n-1$ linearly independent lattice points from the
ball of radius $r=2^{n-1}\oo_{n-1}^{-1}D^{-(n-2)}\rho^{\frac{1}{n-1}}$. From
here $\rho$ can be expressed as a function of $r$,
and the lower bound in Theorem 3 follows.

\demo{Proof of Claim}
We shall assume throughout this argument that $\rho$ is sufficiently large
in terms of $n$. The inequality $\ll_1(v) \leq D\rho^{\frac{1}{n-1}}$
is equivalent to the existence of
a primitive $u \in \ZZ$ with $vu=0$ and $|u| \leq
D\rho^{\frac{1}{n-1}}$. In other words, $v \in\LL(u)$ for
some primitive $u$ with $|u|\leq D\rho^{\frac{1}{n-1}}$.
For any primitive $u$ with $|u|\leq D\rho^{\frac{1}{n-1}}$,
we estimate the number of corresponding vectors $v$.

Using (6) we can see that $\ll_{n-1}(u)\leq
2^{n-1}\oo_{n-1}^{-1}D\rho^{\frac{1}{n-1}}=o(\rho)$ which implies
that $\LL(u)$ contains a lattice parallelotope of nonzero volume
and of diameter $o(\rho)$. Therefore
the number of corresponding vectors $v$ is at most
$$|\LL(u)\cap\rho \BB^n|\leq 2\text{vol}(\rho \BB^{n-1})/\det\LL(u)=2\oo_{n-1}\rho^{n-1}/|u|.$$
Hence the total number of $v \in V$ with $\ll_1(v) \leq
D\rho^{\frac{1}{n-1}}$ is at most
$$2\oo_{n-1}\rho^{n-1}\sum_{|u| \leq D\rho^{\frac {1}{n-1}}}
\frac{1}{|u|} \leq 4\oo_{n-1}\oo_n D^{n-1}\rho^n,
$$
as can be shown by a straightforward calculation. The total
number of points in $V$ is at least
$\frac{1}{2\zeta(n)}\oo_n\rho^n$. Thus, the
number of $v \in V$ with $\ll_1(v) \geq
D\rho^{\frac{1}{n-1}}$ is at least
$$
\left(\frac {1}{2\zeta(n)}-4\oo_{n-1}D^{n-1}\right)\oo_n\rho^n,
$$
which is larger than $\oo_n\rho^n/10$ if the constant $D$ is chosen
properly. \qed
\enddemo

\demo{Proof of Corollary} We have
$$\align
\sum_{H\in\HH_r}|H\cap r\BB^n\cap\ZZ|
&=|\HH_r|+\sum_{0\ne v\in r\BB^n\cap\ZZ}|\{H\in\HH_r|v\in H\}|\\
&\leq|\HH_r|+|r\BB^n\cap\ZZ|^{n-1}\\
&\leq|\HH_r|+\omega_n^{n-1}(r+\sqrt n)^{n(n-1)},
\endalign$$
where the first inequality follows from the fact that $H\in\HH_r$ is
spanned by $v$ and other $n-2$ independent vectors in $r\BB^n\cap\ZZ$.
By Theorem 3 we have
$$s_r=\frac{1}{|\HH_r|}\sum_{H\in\HH_r}|H\cap
r\BB^n\cap\ZZ|\le1+c_1^{-1}\omega_n^{n-1}(1+\sqrt n/r)^{n(n-1)},$$
where the right hand side is bounded by a function of $n$ as required.\qed
\enddemo

\head 5. Epilogue \endhead

Hal\'asz' question studied in this paper is related to the
following famous problem of Littlewood and Offord \cite{9}.
{\it Given $k$ not necessarily distinct complex numbers,
$v_1, v_2,\ldots,v_k$, whose absolute values are at least $1$,
at most how many of the $2^k$ subset sums
$\sum_{i\in I}v_i, \ \ I\subseteq\{1,2,\ldots,k\}$ can
belong to the same open ball of unit diameter?}

Erd\H os \cite{3} proved that for reals the best possible
upper bound was ${k\choose \lfloor k/2\rfloor}$. G. O. H.
Katona \cite{6} and D. Kleitman \cite{7} independently settled
the original question by showing that the same bound is valid
for complex numbers. Shortly after, Kleitman \cite{8}
managed to generalize this theorem to systems of vectors of
absolute value at least $1$ in any Euclidean space $\RR$.
In all cases, the upper bound is attained when all vectors
(numbers) coincide.

Erd\H os and Moser considered the similar problem of haw many subset sums of
$k$ {\it distinct} numbers can coincide. A. S\'ark\"ozy, E. Szemer\'edi
\cite{12} found the order of the magnitude of this number and later
R. Stanley \cite{13} found the exact answer. G. Hal\'asz \cite{5}
considered the similar problem of how many subset sums can coincide under
various assumptions assuring that the $k$ vectors are quite different.
J. Griggs and G. Rote
\cite{4} investigated the following problem of this type.
{\it Given $k$ $n$-wise linearly independent vectors
$v_1, v_2,\ldots,v_k\in\RR$, at most how many of the $2^k$
subset sums $\sum_{i\in I}v_i, \ \ I\subseteq\{1,2,\ldots,k\}$
can coincide?} Denoting this function by $f_n(k),$ they
obtained that
$$f_n(k)>C_n\frac{2^k}{k^{3n/2-1}},$$
and it is implicit in Hal\'asz \cite{5} that
$$f_n(k)<C'_n\frac{2^k}{k^{n/2+\lfloor{n/2}\rfloor}}.$$
(Here $C_n$ and $C'_n$ are positive constants depending only
on the dimension $n$.) The orders of magnitude of these two
bounds differ already in $3$-space ($n=3$).

Note that the construction of Griggs and Rote \cite{4} can
be regarded as the special case of our construction at
the end of Section 2, when $\CC$ is a box of the form
$[0,1]\times[0,x]^{n-1}$.

Hal\'asz observed that the construction in \cite{4} can be
extended to give the following result. Let $h_n(r)$ denote
the maximum number of $n$-wise linearly independent lattice
points that can be chosen in $r\BB^n$. Let $r(k)$ be the
smallest $r$ for which $h_n(r)\geq k$. Then
$$f_n(k)>C''_n\frac{2^k}{k^{n/2}r^n(k)}.$$
This would improve on the previous lower bound, provided
that $r(k)=o(k^{(n-1)/n})$, or, equivalently,
$$\lim_{r\rightarrow\infty}\frac{h_n(r)}{r^{n/(n-1)}}=\infty.$$
However, the results in this paper show that this is not the
case.

With the exception of Erd\H os, all Hungarian mathematicians
mentioned in this section (G\'abor Hal\'asz, Gyula Katona,
Andr\'as S\'ark\"ozy, Endre Szemer\'edi) recently have
turned or will turn {\it sixty}. We congratulate them with
this note.

\enddocument